\documentclass[ijoc,nonblindrev]{informs3} 
\usepackage{savesym}

\DoubleSpacedXII
\usepackage{geometry}
\usepackage[parfill]{parskip}    
\usepackage{graphicx}
\usepackage{changes}
\usepackage{color,soul}
\usepackage{amsmath,amssymb}
\savesymbol{0}
\usepackage{outlines}
\usepackage{lineno}
\usepackage{bbm}
\newcommand{\conv}{\text{conv}}
\newcommand{\inter}{\text{int}}
\newcommand{\bdef}{\stackrel{\triangle}{=}}
\newcommand{\pf}{\textbf{Proof: }}
\newcommand{\R}{\mathbb{R}}
\newcommand{\X}{\mathbb{X}}
\newcommand{\x}{\mathbf{x}}
\newcommand{\z}{\mathbf{z}}
\newcommand{\bmu}{\mathbf{\mu}}
\newcommand{\bnu}{\mathbf{\nu}}
\newcommand{\ie}{{i.e.}}
\newcommand{\eg}{{e.g.}}
\newcommand{\cf}{{cf.\ }}
\newcommand{\one}{\mathbf{1}}
\newcommand{\zero}{\mathbf{0}}
\newcommand{\la}{{\langle}}
\newcommand{\ra}{{\rangle}}
\newcommand{\bb}{{\mathbf{b}}}

\makeatletter
\newcommand{\rmnum}[1] {\romannumeral #1}
\newcommand{\Rmnum}[1]{\expandafter\@slowromancap\romannumeral #1@}
\makeatother

\makeatletter
\def\munderbar#1{\underline{\sbox\tw@{$#1$}\dp\tw@\z@\box\tw@}}
\makeatother

\usepackage{graphicx}
\usepackage{amssymb}

\usepackage{lineno}

\begin{document}



\TITLE{Variable Partitioning for Distributed Optimization}



\ARTICLEAUTHORS{%
\AUTHOR{Yuchen Zheng, \EMAIL{richardzyc@gatech.edu}}
\AFF{H. Milton Stewart School of Industrial \& Systems Engineering, Georgia Institute of Technology}
\AUTHOR{Ilbin Lee, Ph.D., \EMAIL{ilbin@ualberta.ca}}
\AFF{Alberta School of Business, University of Alberta}
\AUTHOR{Nicoleta Serban, Ph.D., \EMAIL{nserban@isye.gatech.edu}}
\AFF{H. Milton Stewart School of Industrial \& Systems Engineering, Georgia Institute of Technology}
} 

\ABSTRACT{%
This paper is about how to partition decision variables while decomposing a large-scale optimization problem for the best performance of distributed solution methods. Solving a large-scale optimization problem sequentially can be computationally challenging. One classic approach is to decompose the problem into smaller sub-problems and solve them in a distributed fashion. However, there is little discussion in the literature on which variables should be grouped together to form the sub-problems, especially when the optimization formulation involves complex constraints. We focus on one of the most popular distributed approaches, dual decomposition and distributed sub-gradient methods. Based on a theoretical guarantee on its convergence rate, we explain that a partition of variables can critically affect the speed of convergence and highlight the importance of the number of dualized constraints. Then, we introduce a novel approach to find a partition that reduces the number of dualized constraints by utilizing a community detection algorithm from physics literature. Roughly speaking, the proposed method groups decision variables that appear together in constraints and solves the resulting sub-problems with blocks of variables in parallel. Empirical experiments on a real application show that the proposed method significantly accelerates the convergence of the distributed sub-gradient method. The advantage of our approach becomes more significant as the size of the problem increases and each constraint involves more variables.  
}%


\KEYWORDS{Large-scale optimization; Decomposition; Distributed computing; Community structure detection }

\maketitle



\section{Introduction}\label{sec:intro}
Solving large-scale optimization problems using one computer core and sequential computing can be computationally challenging due to the data storage and retrieval, and due to the computational load and memory usage for obtaining an optimal solution. Distributed computing is a popular framework for tackling the computational complexity of large-scale optimization \cite{androulakis_et_al1996,boyd_et_al2011,camponogara_oliveira2009,inalhan_et_al2002,nedic_ozdaglar2009a,palomar_chiang2006,raffard_et_al2004,richtarik_takac2016,shastri_et_al2011,simonetto_jamalirad2016,terelius_et_al2011,xiao_et_al2004}. Distributed computing for optimization problems involves two computational considerations: the decomposition into smaller sub-problems in a way that each sub-problem can be stored and solved in a single machine, and the derivation of a solution for the original optimization problem by (iteratively) solving the sub-problems \cite{boyd_et_al2011,camponogara_oliveira2009,nedic_ozdaglar2009a,palomar_chiang2006,raffard_et_al2004,shastri_et_al2011,simonetto_jamalirad2016,terelius_et_al2011,xiao_et_al2004}. Applications of distributed optimization arise in various emerging areas, such as resource allocation over large-scale networks \cite{nedic_ozdaglar2009a,palomar_chiang2006,xiao_et_al2004}, aircraft coordination \cite{inalhan_et_al2002,raffard_et_al2004}, and estimation problem in sensor networks \cite{duchi_et_al2012}.

One classic approach to the decomposition of an optimization problem is the {\it dual decomposition}. It decomposes an optimization problem into smaller sub-problems by relaxing some of the constraints. Then, the resulting Lagrangian dual is solved by a distributed sub-gradient algorithm \cite{palomar_chiang2006,raffard_et_al2004,simonetto_jamalirad2016,terelius_et_al2011,xiao_et_al2004}. Another decomposition technique is introducing copy variables and the alternating direction method of multipliers (ADMM) is commonly used in the distributed algorithm \cite{boyd_et_al2011}. Most of the existing works have focused on developing either a new decomposition technique or a novel distributed algorithm. However, regardless of a decomposition technique (e.g., dual decomposition) or a distributed solution algorithm (e.g., sub-gradient method), partitioning the decision variables across sub-problems remains one of the key challenges in distributed optimization. For instance, for a network optimization problem over a graph, should we define a sub-problem for each node or a group of nodes? If a sub-problem may contain multiple nodes, how should we assign nodes to sub-problems? Past works discussing this issue are limited. In \cite{parikh2014block}, a distributed block splitting algorithm based on graph projection splitting was introduced for decomposing and solving large-scale problems in parallel in which the objective function is separable by blocks of variables. However, the authors pointed out that, in practice, it was not obvious which subset of variables should be processed together rather than on separate machines. Thus, the key question of this paper still remained unresolved. There have been other efforts to determine how a complex system should be partitioned to achieve faster convergence, spectral clustering technique in \cite{guo2016intelligent} and simultaneous partitioning and coordination strategy in \cite{allison2009optimal}. However, these works were focused on specific applications with relatively small numbers of variables (a few hundreds or less).

The computational approach in this paper is motivated by the observation that a decomposition of an optimization problem (in other words, a partition of decision variables) critically affects the computational performance of distributed optimization algorithms. For illustration, we used one of the most common approaches in distributed optimization, dual decomposition and distributed sub-gradient method. Sub-gradient methods have been shown to converge to optimality as long as the resulting Lagrangian dual satisfies strong duality,  regardless of which constraints are dualized or how decision variables are partitioned \cite{bertsekas1999}. However, our empirical analysis shows that the convergence may be extremely slow, potentially not reaching convergence even after a large number of iterations (e.g., ten thousands), especially when there is a large number of highly complex constraints. 

In this paper, we provide a theoretical explanation of why a partition of decision variables can affect the convergence of distributed sub-gradient methods. In relation to the theoretical upper bound for the convergence rate of sub-gradient methods established in literature \cite{boyd_note,goffin1977}, we explain the importance of minimizing the number of dualized constraints to achieve faster convergence of distributed sub-gradient methods. The intuition is that the more ``similar'' the Lagrangian dual and the original problems are, the more desirable it is for the empirical performance of sub-gradient methods, yielding faster convergence results.

A key contribution of this paper is a novel method to find a partition of decision variables for dual decomposition and to solve the resulting sub-problems with blocks of variables. Our focus is developing a method that dualizes as fewer number of constraints as possible while decomposing a large-scale optimization problem. Our approach to find such a partition uses a \emph{community detection} algorithm from the physics literature \cite{newman2004,newman2006}. The goal of community detection is to identify community structures within a network, in other words, to find groups of nodes in such a way the connections within each group are dense while there is little connectivity between the groups. It has been applied to the Internet, citation networks, social networks among others \cite{fortunato2010community}. Roughly speaking, we use community detection to group decision variables that tend to appear in constraints together so that the number of constraints that involve variables over multiple sub-problems (thus, need to be dualized) is minimized.

The proposed approach is general and applicable to various problem classes, but for illustration purpose, we present the proposed method applied to transportation problems as follows. First, we construct a graph whose nodes represent demand locations. Two nodes are connected if there is a constraint involving the two demand locations (e.g., a supply location can serve the two demand locations and there is a capacity limit constraint at the supply location). Each edge is weighted by the number of constraints involving the two demand locations. Then, we apply a community detection algorithm to the graph to find a partition of demand locations. The community detection algorithm identifies communities of demand locations such that demand locations in the same community are densely `connected' by the constraints but those in different communities are sparsely `connected' by the constraints. Then, we decompose the original optimization problem into blocks of demand locations given by the community detection, thus reducing the number of dualized constraints by utilizing the community structure. Our empirical illustration in Section~\ref{sec:results} shows that the method introduced in this paper significantly accelerates the convergence of distributed sub-gradient methods. 

This paper is organized as follows. First, we review dual decomposition and sub-gradient methods in Section~\ref{sec:background}. In this section, we also analyze why decomposition is important for the performance of sub-gradient methods. In Section~\ref{sec:main}, we introduce our approach to find a decomposition via community detection. We illustrate its empirical performance for a real application in Section~\ref{sec:results} and conclude in Section~\ref{sec:conclusion}. 
\section{Dual Decomposition and Sub-gradient Method}\label{sec:background}
Dual decomposition is a common technique for decomposing a large-scale optimization problem into smaller sub-problems \cite{bertsekas1999,raffard_et_al2004,terelius_et_al2011}. Given a partition of decision variables, constraints that are over multiple groups of variables are relaxed and added to the objective function as penalty terms for violation, so that the Lagrangian relaxation is decomposable into smaller sub-problems. In this section, we first review the dual decomposition technique, followed by a distributed sub-gradient algorithm. We also analyze its convergence rate established in the existing literature and discuss why the convergence may be slow.

\subsection{Transportation Problem}

The transportation problem is a general class of problems, in which commodities are transported from a set of sources to a set of destinations. Let $x_{ij}$ denote the matching variable from demand location $i\in I$ to supply location $j\in J$. Let $X$ denote the $|I|\times |J|$ matrix whose $(i,j)$ entry is $x_{ij}$ and $X_i$ denotes the $i$th row. The general  optimization model is given as follows.
\begin{align}
\text{(GP)}\ \min_{X}\ \sum_{i \in I} &\sum_{j \in J}w_{ij} x_{ij}\label{eq:GeneralModelObjective}\\
\text{s.t. }
\sum_{j \in J_i} x_{ij} &\ge m_i \text{ for }i\in I,\label{eq:GeneralModelLocalConstraintI}\\
\sum_{i \in I_j} x_{ij} &\le s_j \text{ for }j\in J,\label{eq:GeneralModelLocalConstraintJ}\\
X&\ge 0, \label{eq:NonNegativeConstraints2}
\end{align} 
where $m_i$ is the minimum demand that needs to be satisfied at each demand location $i \in I$, $s_i$ is the maximum capacity at each supply location $j \in J$, $w_{ij}$ is the cost associated with demand location $i$ getting one unit of goods from supply location $j$, $J_i$ is the set of supply locations that can serve demand location $i$, and $I_j$ is the set of demand locations that can be served by supply location $j$. 
In real applications where there is a large number of demand and supply locations, it is often assumed that each demand location can only be served by a subset of supply locations. For instance, in logistics, suppliers may have access only to a few demand locations due to regions of operations, or that some demand locations are simply too far away. 
In this paper, we consider only continuous decision variables. For example, $x_{ij}$ may be a number of service hours assigned to demand location $i$ from supply location $j$.  

\subsection{Dual Decomposition and Distributed Sub-gradient Method}\label{sec:baseline}

 We first review dual decomposition of (GP) with a straightforward partition of decision variables, a sub-problem for each demand location $i$. In (GP), constraints \eqref{eq:GeneralModelLocalConstraintJ} include variables over multiple demand locations. In order for (GP) to be decomposed for each $i$, the constraints \eqref{eq:GeneralModelLocalConstraintJ} are relaxed and appended as penalties for their violation to the objective function. Let $\lambda_j \ge 0$ be the dual variable for each constraint in \eqref{eq:GeneralModelLocalConstraintJ}. The resulting Lagrangian is:
\begin{align*}
L(X,  \Lambda)=&\sum_{i \in I} \sum_{j \in J_i}w_{ij} x_{ij}+\sum_{j\in J}\lambda_j (\sum_{i\in I_j}x_{ij}-s_j)=\sum_{i\in I}\sum_{j\in J_i}(w_{ij}x_{ij}+\lambda_j x_{ij}) -\sum_{j \in J}\lambda_j s_j\\
=&\sum_{i\in I}L_i(X_i,  \Lambda)-\sum_{j \in J}\lambda_j s_j,
\end{align*}
where
\begin{equation*}
L_i(X_i,  \Lambda)\triangleq \sum_{j\in J_i}(w_{ij}x_{ij}+\lambda_j x_{ij}).
\end{equation*}
Given that $\lambda_j$ are fixed, the Lagrangian is decomposed for each demand location $i$.

Let $\mathcal{D}$ be the set of values of $X$ satisfying the constraints of (GP) that are not dualized, \ie, \eqref{eq:GeneralModelLocalConstraintI} and \eqref{eq:NonNegativeConstraints2}. Let $\mathcal{D}_i$ be the set of values of $X_i$ that satisfy those constraints restricted to demand location $i$. 
Then, the dual objective function $g(\Lambda)$ can be computed by optimizing the decomposed Lagrangian for each demand location separately as follows:
\begin{equation*}
g(\Lambda)\triangleq \inf_{X\in \mathcal{D}}L(X, \Lambda)=\sum_{i\in I}g_i(\Lambda),
\end{equation*}
where
\begin{equation*}
g_i(\Lambda)\triangleq \inf_{X_i\in \mathcal{D}_i}L_i(X_i,\Lambda).
\end{equation*}
The local optimization problem on the right hand side can be written as follows:
\begin{align*}
\text{(LR$_i$)}\ \min_{X_i}\ L_i(X_i,\Lambda)
&=\sum_{j\in J_i}(w_{ij}x_{ij}+\lambda_j x_{ij})\\
\text{s.t. }\sum_{j\in J_i} x_{ij} &\ge m_i,\\
X_i&\ge 0.
\end{align*}

The Lagrangial dual of (GP) is defined as
\begin{equation}
\text{(LD)}\ \sup_{\Lambda\ge 0}g(\Lambda).
\end{equation}

The optimal values of (GP) and (LD) coincide (\ie, strong duality holds) because the objective function of (GP) is convex and the constraints of (GP) are affine, thus, it satisfies Slater's constraint qualification \cite{boyd2004convex}. 

A distributed sub-gradient algorithm for solving (LD) is given as follows.

 \begin{center}
 \textbf{Distributed Sub-gradient Algorithm}
 \end{center}
 \begin{outline}
 \1[1.] Choose a starting point $\Lambda^1$. Let $t:=1$ (first iteration).
 \1[2.] Solve the local optimization problem (LR$_i$) with $\Lambda=\Lambda^t$ for each demand location $i\in I$ to obtain $X_i^t$. 
 \1[3.] If a given stopping criterion is satisfied, stop. Otherwise, $t:=t+1$, update the multipliers as below, and go to Step 2:
 \begin{align}
 \lambda^{t+1}_j&=\max\{\lambda^t_j+\alpha_t (\sum_{i\in I_j}x^t_{ij}-s_j), 0\}\text{ for }j\in J.
 \end{align}
 \end{outline}

It is well-known that if the step-size $\{\alpha_t\}^\infty_{t=1}$ satisfies
\begin{equation}\label{eq:stepsize}
\sum_{t=1}^\infty \alpha_t =\infty\text{ and }\sum_{t=1}^\infty \alpha_t^2 <\infty,
\end{equation}
then the value of $g(\Lambda^t)$ converges to the optimal objective function value of (GP) (\eg, see \cite{bertsekas1999}). Moreover, the running average of the primal iterates $X^t$ becomes optimal for (GP) asymptotically as $t$ goes to infinity \cite{nedic_ozdaglar2009b,simonetto_jamalirad2016}.

\subsection{Analyzing Convergence Rate of Sub-gradient Method}\label{sec:rate_analysis}

Convergence rates of sub-gradient methods have been established under various settings \cite{boyd_note,goffin1977}. We first review a convergence rate result of the sub-gradient method with the step-size given in \eqref{eq:stepsize} and discuss its slow convergence and our proposed approach to address it. 

\begin{figure}
	\centering
	\includegraphics[width=\textwidth]{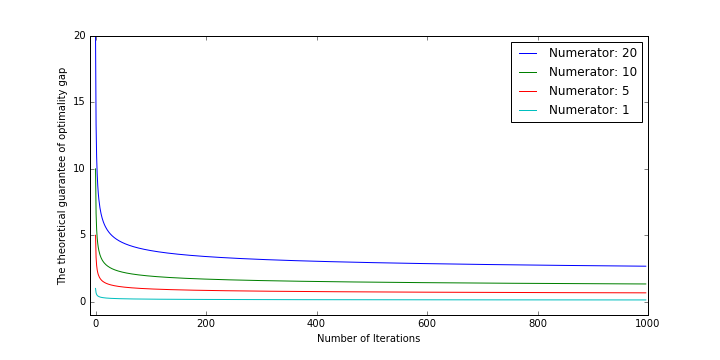}
	\caption{Convergence of the theoretical guarantee of optimality gap with different numerator values and $\alpha_t=\frac{1}{t}$.}
	\label{convGap}
\end{figure}

Since sub-gradient methods do not improve monotonically, it is common to keep track of the best solution up to the current iteration. Let $\Lambda^t_{\text{best}}$ denote the solution having the lowest $g$ value at the end of iteration $t$. Let $R$ be an upper bound on the distance between the initial dual solution and the set of optimal dual solutions, \ie, $||\Lambda^1-\Lambda^\star||_2\le R$. Also, let $G$ be an upper bound on the norm of the sub-gradients computed by the algorithm, \ie, $||h^t||_2\le G$, where $h^t\in \R^{|J|}$ and $h^t_j=\sum_{i\in I_j}x^t_{ij}-s_j$ for $j\in J$. From \cite{boyd_note}, we have the following upper bound on the optimality gap, which goes to zero as $t$ goes to infinity: 
\begin{equation}\label{eq:convergence_rate}
g(\Lambda^t_{\text{best}})-g(\Lambda^\star) \le \frac{R^2+G^2||\alpha||_2^2}{2\sum_{k=1}^t\alpha_k}.
\end{equation}
However, depending on the value of its numerator, the upper bound may converge to zero so slowly that it does not approach zero even at a large value of $t$ (e.g., hundreds of thousands). See Figure \ref{convGap} illustrating the significant difference in the convergence of the upper bound depending on the value of the numerator where $\alpha_t=\frac{1}{t}$. More importantly, our experimental results for a transportation problem (Section~\ref{sec:results}) show that the optimality gap itself may not approach zero even after a large number of iterations under the dual decomposition for each demand location. We emphasize that the slow convergence of the theoretical upper bound applies to any optimization problem, not limited to transportation problems used for illustration in this paper.

We will next investigate possible ways to speed up the convergence of the upper bound. The upper bound contains the step size $\{\alpha_t\}$, an upper bound $R$ on the distance between the initial dual solution and the optimal dual solution set, and an upper bound $G$ on the magnitude of the sub-gradients. Adjusting the step size is an easy choice for accelerating sub-gradient methods, but it is purely empirical and specific to each application. Another important component that governs the behavior of the upper bound is \emph{the number of dualized constraints}, because it equals the dimension of dual parameter vector $\Lambda$ and the dimension of the sub-gradients. Therefore, the number of dualized constraints is closely related to the magnitude of $R$ and $G$, and thus, directly affects the convergence of the upper bound. 

In addition, note that each component of the sub-gradient at a primal solution $X$ is the violation of the corresponding dualized constraint at $X$. Thus, the fewer violations of the dualized constraints are, the lower the magnitude of the sub-gradient is. This again emphasizes the importance of the number of dualized constraints. Dualizing more constraints leads to a more relaxed feasible region of the resulting Lagrangian relaxation. Then, the primal iterates obtained while running the sub-gradient method have more ``room'' to deviate from the original feasible region, thus allowing larger violations of dualized constraints, which leads to higher magnitudes of the sub-gradients, and thus, a higher $G$ value and slower convergence.

For the transportation problem and the dual decomposition introduced in the previous section, the above intuition is translated as follows. 
The Lagrangian relaxation of (GP) was obtained by dualizing all of the supply constraints \eqref{eq:GeneralModelLocalConstraintJ}, so the resulting relaxation differs significantly from the original problem. The subproblem (LR$_i$) for demand location $i$ is simply matching the demand of $i$ to accessible supply locations where the values of the dual multipliers make the location $i$ prefer some supply locations than others. Thus, the competition among demand locations for limited resources is only indirectly reflected via the dual multipliers. In other words, the level of decomposition is so fine that each sub-problem (LR$_i$) loses an important aspect of the original problem, which makes the overall convergence slow. 

A partition of decision variables critically affects the computational performance of the sub-gradient method, and the goal of this paper is to develop a novel method to find a decomposition that speeds up distributed algorithms. For this purpose, we aim at dualizing as fewer number of constraints as possible while taking the computational advantage of distributed computing by dual decomposition. Consider decomposing (GP) into a certain number of sub-problems by partitioning demand locations. The demand constraints \eqref{eq:GeneralModelLocalConstraintI} are decomposable by demand locations, but the supply constraints \eqref{eq:GeneralModelLocalConstraintJ} are not. Given a partition of demand locations, those supply constraints involving demand locations in multiple groups need to be dualized in order for the remaining constraints to be decomposable. Herein we define two demand locations \emph{connected} if and only if there exists a supply location that can serve both of the demand locations, \ie, they appear together in the supply constraint. Then, finding a decomposition with a minimal number of dualized constraints translates into finding a partition in which demand locations in the same group are closely connected and those from different groups are loosely connected. We expand on this approach in the next section. 


\section{Community Detection and Block Dual Decomposition}\label{sec:main}
In this section, we introduce a novel framework for dual decomposition and distributed optimization. The framework consists of community detection and block dual decomposition. Our approach utilizes the structure of decision variables and constraints in order to speed up the convergence of distributed sub-gradient methods. We illustrate the general idea using the transportation problem, but we also introduce a more general version of our approach in Section~\ref{sec:general}. 

\subsection{Community Detection}\label{sec:community}

In this sub-section, we explain how we use community detection to find a partition of decision variables and review the community detection algorithm we used in this paper.

We first build a network graph of demand locations. 
Consider a graph of $n$ nodes, with each node representing one demand location. Two nodes are connected by an edge if the corresponding demand locations are connected, that is, the two demand locations have access to a common supplier. The edge is weighted by the number of suppliers that can serve both of the locations, \ie, the number of constraints the two demand locations appear together. To this network, we apply a community detection algorithm to identify communities of demand locations where those in the same community are densely connected and those in different ones are sparsely connected. Then, we decompose the optimization problem according to the communities. 

Among various algorithms developed in the community detection literature, we use the fast hierarchical agglomeration algorithm proposed by Clauset, Newman, and Moore \cite{clauset2004finding}. The computational complexity of the algorithm is linear in the size of the network for many real-world networks. We briefly explain how the algorithm works below. 

The community detection algorithm is based on a measure of a partition called the {\it modularity}, which evaluates how dense connections are within communities and how few there are between communities \cite{newman2004finding}. Before defining the measure, we introduce some notation. An $n$-by-$n$ weighted adjacency matrix $C$ is defined as
\begin{equation*}
C_{vw}=\begin{cases}
e_{vw} \text{	if nodes $v$ and $w$ are connected,}\\
0 \text{	otherwise,}\\
\end{cases}
\end{equation*}
where $e_{vw}$ is the weight of the edge $(v, w)$. Consider a partition of the nodes and for a node $v$, let  $c_v$ denote the community to which $v$ belongs. Let $\delta(c_v, c_w)$ be 1 if $c_v=c_w$ and 0 otherwise. Let $m=\frac{1}{2}\sum_{v,w}C_{vw}$ be the sum of weights of all edges in the graph and let $k_v=\sum_w C_{vw}$ be the sum of weights of all edges from $v$. Then, the modularity of a partition is defined as:
\begin{equation}\label{eq:modularity}
Q=\frac{1}{2m} \sum_{v,w}\left(C_{vw}-\frac{k_v k_w}{2m}\right)\delta(c_v, c_w).
\end{equation}
An interpretation of the modularity $Q$ is in order. The fraction $k_v k_w/(2m)$ is the expected number of edges between $v$ and $w$ where $m$ edges are randomly assigned between nodes. Thus the modularity measures how strong the community structure is over a random assignment of edges. 
More details of the modularity measure can be found in \cite{clauset2004finding,newman2004finding,newman2004}. In practice, networks with the modularity greater than 0.3 appear to indicate significant community structure \cite{newman2004}.

The community detection algorithm starts with a trivial division where each of the demand location forms a community. Then, it repeatedly joins two communities that results in the biggest increase of the modularity, until it reaches a partition where none of the join operations improves the modularity score. More details of the algorithm can be found in \cite{clauset2004finding}. 

\subsection{Block Dual Decomposition}
Each community in the output of the community detection may involve multiple demand locations. Thus, the corresponding dual decomposition yields sub-problems that include blocks of demand locations, and thus we call the proposed approach \emph{block dual decomposition}. However, we emphasize that the blocks are not arbitrarily created, but we grouped those that are closely connected. This characteristic of our approach makes the resulting sub-problems keep as much structure of the original optimization problem as possible, which is critical for the performance of the distributed sub-gradient method as explained in Section~\ref{sec:rate_analysis} and empirically shown in Section~\ref{sec:results}.

Let $I_b$ for $b=1,\ldots,B$ be the partition of demand locations given by the community detection, thus satisfying $\cup_{b=1}^B I_b=I$ and $I_b\cap I_{b'}=\emptyset$ for $b\neq b'$. Let $J_b$ be the set of supply locations that can serve the demand locations in $I_b$ (\eg, within a pre-specified distance). Note that the set $J_b$'s may not be disjoint as opposed to $I_b$'s. For each block $b$, the suppliers that can serve only the demand locations in the block are said to be {\it interior suppliers} of block $b$, denoted as $J_b^{\text{in}}$, and the suppliers that are not interior suppliers but can serve a demand location in $I_b$ are called {\it boundary suppliers} of block $b$, denoted as $J_b^{\text{out}}$. Let $J^{\text{in}}=\cup_{b=1}^B J_b^{\text{in}}$ and $J^{\text{out}}=\cup_{b=1}^B J_b^{\text{out}}$. Note that $J_b^{\text{in}}$ for $b=1,\ldots,B$ are disjoint. For a demand location $i$, let $b(i)$ denote the block to which $i$ belongs.

Consider the following Lagrangian relaxation of (GP):
\begin{align}
\text{(BLR)}\ \min_{X_i}\ L(X,  \Lambda)=&\sum_{i \in I} \sum_{j \in J_i}w_{ij} x_{ij}+\sum_{j\in J^{out}}\lambda_j (\sum_{i\in I_j}x_{ij}-s_j)\label{eq:BLRObjective}\\
\text{s.t. }\sum_{j\in J_i} x_{ij} &\ge m_i \text{ for }i\in I,\label{eq:BLRLocalConstI}\\
\sum_{i\in I_j} x_{ij} &\le s_j \text{ for }j\in J^{in},\label{eq:BLRLocalConstJ}\\
X&\ge 0.
\end{align}
Note that among the supply side constraints \eqref{eq:GeneralModelLocalConstraintJ}, only those corresponding to boundary suppliers were dualized in (BLR). Consequently, a fewer number of dual variables are needed than in the previous section. By following similar steps to those of the previous section, (BLR) is decomposed as follows:
\begin{align*}
\text{(BLR$_b$)}\ \min_{X_i}\ \sum_{i\in I_b} \sum_{j \in J_i}&w_{ij} x_{ij}+\sum_{i\in I_b}\sum_{j\in J_b^{out}\cap J_i}\lambda_j x_{ij}\\
\text{s.t. }\sum_{j\in J_i} x_{ij} &\ge m_i \text{ for }i\in I_b,\\
\sum_{i\in I_j} x_{ij} &\le s_j \text{ for }j\in J_b^{in},\\
X_i&\ge 0 \text{ for }i\in I_b.
\end{align*}

The resulting distributed subgradient algorithm is as follows.
 \begin{center}
 \textbf{Distributed Subgradient Algorithm with Block Dual Decomposition}
 \end{center}
 \begin{outline}
 \1[1.] Choose a starting point: $\Lambda^1=\zero$. Let $t:=1$.
 \1[2.] Solve the local optimization problem (BLR$_b$) for each demand block $b$ to obtain $X_i^t$ for $i\in I_b$. 
 \1[3.] If (some stopping criterion) is satisfied, stop. Otherwise, $t:=t+1$, update the multipliers as below, and go to Step 2:
 \begin{align}
 \lambda^{t+1}_j&=\max\{\lambda^t_j+\alpha_t \left(\sum_{i\in I_j}x_{ij}-s_j\right), 0\}\text{ for }j\in J^{out}.
 \end{align}
 \end{outline}

\subsection{A General Approach}\label{sec:general}

We have illustrated details of the block dual decomposition with community detection under the transportation problem setting. In this section, we present a similar approach, but with a broader applicability. Consider the following convex optimization problem:
\begin{align*}
\min_{\x}&\  f(\x)\\
\text{s.t. }& A \x \le b, \\
\end{align*}
where $A \in R^{m\times n}$ and $f: R^n \rightarrow R$ is decomposable for each component of $x$, i.e., $f(\x)=\sum_{i =1,...,n} f_i(x_i)$. 
For this general formulation, we illustrate how our approach can be applied to find a partition of decision variables for which the corresponding dual decomposition dualizes a minimal number of constraints.

Construct a graph in which each node represents a decision variable. Two nodes are connected if the two decision variables appear together in a constraint and the edge is weighted by the number of constraints they appear together. Note that in this section, each node represents a decision variable as opposed to the previous section where each node corresponds to a demand location for the transportation problem. Then, we apply the community detection algorithm to this network in order to identify a partition of decision variables where connections within a group are dense but those between groups are sparse. A weighted adjacency matrix $C$ is constructed as follows. We first form an indicator matrix $\tilde A \in R^{m\times n}$ such that for all $i=1,...,m$ and $j=1,...,n$,
\begin{equation*}
\tilde A_{ij}=\begin{cases}
1 \text{	if $A_{ij}> 0$ or $A_{ij}<0$}\\
0 \text{	if $A_{ij}=0$}.\\
\end{cases}
\end{equation*}
Thus, $\tilde A_{ij}=1$ if $x_i$ appears in constraint $j$. Then, an $n$-by-$n$ weighted adjacency matrix $C$ is defined as $C=\tilde{A}\tilde{A}^T$, thus $C_{uv}$ is the number of times variables $x_u$ and $x_v$ appear in the same constraint, for all $u=1,...,n$ and $v=1,...,n$. Then, the modularity score of a partition is computed by using this $C$ matrix as \eqref{eq:modularity} and the community detection algorithm is applied.
If the objective function is decomposable by groups of decision variables instead of each individual variable, then the aforementioned algorithm can be trivially extended by treating each block of variables as one node in the graph.

%

\section{Numerical Results}\label{sec:results}
In this section, we present experimental results for our approach with application to the transportation problem. We first explain the application and problem generation setup. Then, we empirically compare the sub-gradient method with the dual decomposition for each demand location (Section~\ref{sec:background}) and our approach (Section~\ref{sec:main}) for problem instances with varying sizes and network structures. The solver was implemented in Julia, a high-performance dynamic programming language for numerical computing \cite{bezanson2012julia}, along with Gurobi for optimization.  

\subsection{Problem Setup}
For experiments, we generated problem instances based on an optimization model from a real application: matching children in need of healthcare to care providers in Georgia. The optimization problem is in the form of (GP) and it minimizes the total travel distance that the patients in each census tract have to travel to receive care. Each census tract is a demand location and each provider location is a demand location. Providers' practice location addresses, which are supply locations, were obtained from the 2013 National Plan and Provider Enumeration System (NPPES). The patient population is aggregated at the census tract level using the 2010 SF2 100\% census data. In order to compute the number of children in each census tract, the 2012 American Community Survey data was used. More details about the application problem can be found in \cite{gentili_et_al2017}. 

The optimization problem is (GP), where $x_{ij}$ denotes the number of children in demand location $i\in I$ assigned to supply location $j\in J$; $w_{ij}$ is the distance between the demand location $i$ and supply location $j$; $m_i$ is the minimum number of patients needed to be served at demand location $i$; $s_j$ is the maximum number of patients supply location $j$ can accommodate. We allow the variables $x_{ij}$ to be fractional for the computational tractability of the problem, and also because the number of children to be assigned from each location is typically very large (approximately 2500-8000 children). In this application problem of Georgia, there are $I=1955$ demand locations and $J = 3157$ supply locations.

Using the optimization problem from real data, we created problem instances with different sizes. First, we divided the map of Georgia into 50 blocks, 10 horizontally by 5 vertically, based on the longitudinal and latitudinal coordinates. Then, we counted the number of census tracts and provider locations in each block. For each block, we constructed a histogram of the demands ($m_i$'s) of the census tracts in the block. We also obtained a histogram of supply capacities ($s_j$'s) of providers for each block. Then, we generated a problem instance for a given number of demand and supply locations as follows. We determined the number of demand and supply locations in each block in a way that the numbers of locations in different blocks of a generated instance are proportional to the numbers of locations in blocks of the original problem. Positions of demand and supply locations in each block were sampled randomly from the uniform distribution over the block. For each demand or supply location, the amount of demand or capacity was sampled from the empirical histogram of the block for demand or supply, respectively. In addition, a demand location $i$ was said to have access to a supply location $j$ if the distance $w_{ij}$ between them is less than a given threshold $d_{\max}$. By changing the threshold $d_{\max}$ on the traveling distance, we were able to adjust the connectivity between demand and supply locations, thus changing the structure of the network. A lower $d_{\max}$ indicates a sparser network. A dummy supply location was included to guarantee feasibility. All demand locations have access to the dummy but with a very large distance, in our case 1000 miles.

\subsection{Comparative Results}
For the generated problem instances, we compared the empirical performance of the sub-gradient method with the two dual decompositions, the one for each demand location in Section~\ref{sec:baseline} (which we call `baseline' for simplicity) and the proposed block dual decomposition via community detection (called `block'). We first implemented the sub-gradient method with the two decomposition methods in a sequential fashion, that is, all the sub-problems in each iteration are solved sequentially using one computing node. In addition, we implemented a parallel version of the sub-gradient method for the block dual decomposition (called `distributed block') in a distributed computing framework. The parallel implementation solves the sub-problems (BLR$_b$) in parallel at each iteration using 3 computing cores (Intel Core Haswell Processors). The step size $\alpha_t$ was chosen to be $c/t$, where $c$ is a constant scaling factor. 
For all of the methods, we have tried different values of $c=1/10, 1/50, 1/80,$ and $1/100$, but all of the methods had the fastest convergence for the same $c$ value at $1/80$, which we used for this comparison. We measured the number of iterations and the CPU run time in seconds required to reach a certain optimality gap percentage. 

\begin{figure}
	\includegraphics[width=\textwidth]{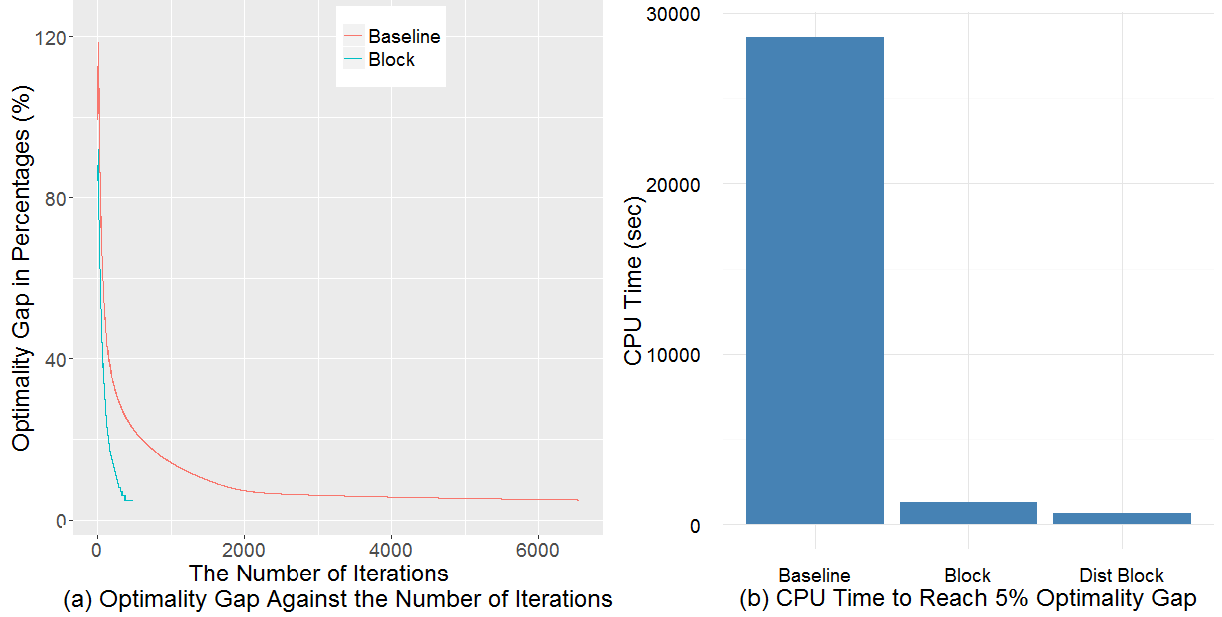}
	\caption{Comparison between the baseline dual decomposition and the block dual decomposition for 1000 demand locations and 1000 supply locations.}
	\label{size}
\end{figure}

Figure~\ref{size} shows the comparison for a problem instance with 1000 demand locations, 1000 supply locations, and $d_{max}=20$(miles), which was generated as previously explained. The instance has 81,230 decision variables in total. Figure~\ref{size}(a) shows how the optimality gap progressed as a function of the number of iterations for the baseline and the block. Both methods were terminated when the dual objective function value was within 5\% of the true optimal objective function value. 

The proposed approach requires significantly fewer iterations to achieve the same optimality gap than the baseline decomposition. The block and the distributed block yield the same performance in the number of iterations to achieve a certain optimality gap of course, but they differ in run time. Figure~\ref{size}(b) shows the CPU run times in seconds to reach 5\% optimality gap for the baseline, the block, and the distributed block. The baseline reached the stopping criterion after 6,654 iterations, about 8 hours. In comparison, the block finished after 496 iterations and 22 minutes. The distributed block finished after 11 minutes, achieving a 44 times speed up comparing to the baseline. 

Table~\ref{table:summary} shows similar comparisons with more problem instances of varying sizes. As the size of the problem grew, 
each algorithm took more iterations and more time to reach 5\% optimality gap and the discrepancies between the methods also grew. The distributed block consistently achieved faster convergence in both the number of iterations and the run time by a large margin. The distributed block yielded 2-2.5 times speed up comparing to its sequential counterpart (the block) using three computing nodes.

\begin{table}[tbp]

	\begin{tabular}{|c|c|c|c|c|c|c|c|c|}
		\hline
		\multicolumn{3}{|c|}{\textbf{Problem Size}}                        & \multicolumn{2}{c|}{\textbf{\# Iterations}} & \multicolumn{3}{c|}{\textbf{Run Time (in Sec.)}} &     \\ \hline
		\textbf{\# Dem} & \textbf{\# Sup} & \textbf{\# Vars}  & \textbf{Baseline}        & \textbf{Block}        & \textbf{Baseline} & \textbf{Block} & \textbf{Dist Block} & \textbf{\# Blocks}\\ \hline
		500                 & 500                 & 21,612                  & 93                   & 47                    & 403           & 183            & 76                  & 21\\ \hline
		1000                & 1000                & 81,230                & 6,654                & 496                   & 28,560        & 1,315          & 655                 & 10 \\ \hline
		1500                & 1500                & 186,062                & 27,190               & 3,009                 & 316,924       & 21,993         & 8,992               & 9\\ \hline
	\end{tabular}

	\caption{Comparison on reaching 5\% optimality gap for problem instances with varying sizes.}\label{table:summary}
\end{table}

\begin{figure}
	\includegraphics[width=\textwidth]{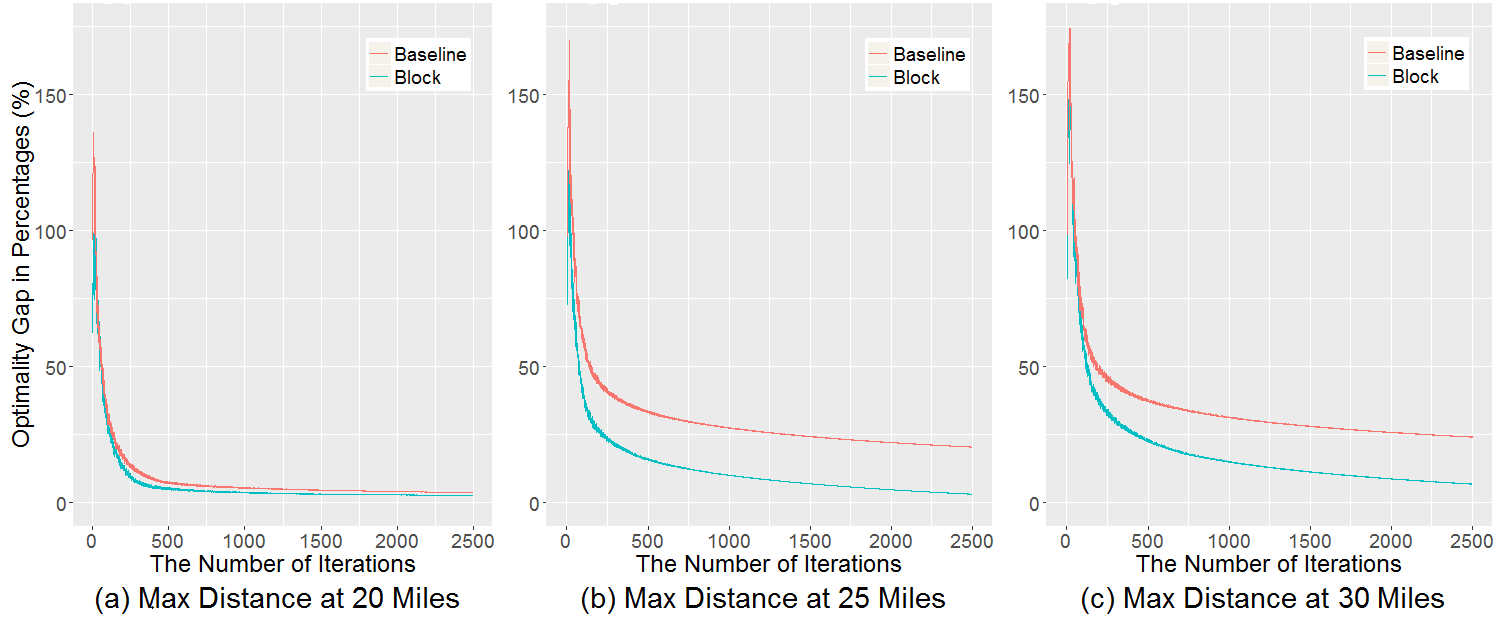}
	\caption{Comparison on rate of convergence between the dual decomposition and the block dual decomposition with varying network structures.}
	\label{sparsity}
\end{figure}

In order to evaluate how the connectivity of locations affects the performance of the methods, we constructed problem instances with 1500 demand locations and 500 supply locations, but different values of $d_{\max}=20, 25$, and $30$ in miles. Recall that $d_{\max}$ is the maximum distance to travel and that the smalller the value is, the more sparse the network is. Figure~\ref{sparsity} compares the baseline and the block for the three instances within 2500 iterations. For the three values of $d_{\max}=20,25$, and $30$, each demand location had access to 40, 57, and 73 providers on average, respectively. 

We observe that the block performs better than the baseline consistently for different values of $d_{\max}$, but the discrepancy of performance gets bigger as the threshold increases, that is, as the network gets more dense. At the $2500$th iteration, the optimiality gap of the baseline is 3.6\%, 20.3\%, and 24.1\% and the optimiality gap of the block is 2.6\%, 3.0\%, and 6.8\% for $d_{\max}=20, 25$, and $30$ miles, respectively. Thus, the advantage of the proposed approach becomes more significant as the network becomes more dense.

The effect of the network structure on the performance can be explained geometrically as follows. For a larger value of the threshold, each provider is accessible from more demand locations and thus, each provider constraint contains more decision variables. In that case, dualizing each provider constraint results in a bigger change on the feasible region in the following sense. For example, imagine the following two relaxations: relaxing $x_1+x_2\le 1$ from $\{(x_1,x_2)\ |\ x_1+x_2\le 1, x_1\ge 0, x_2\ge 0\}$ and relaxing $x_1\le 1$ from $\{(x_1,x_2)\ |\ x_1\le 1, x_1\ge 0, x_2\ge 0\}$. The former yields a bigger change than the latter. Therefore, when $d_{\max}$ is larger, dualizing each provider constraint causes a bigger change on the feasible region. Moreover, note that the baseline dual decomposition dualizes more provider constraints than the proposed approach. Therefore, as $d_{\max}$ increases, the discrepancy between the feasible regions of the Lagrangian relaxation and the original problem becomes more significant for the baseline than it does for the proposed approach. Thus, when $d_{\max}$ increases (\ie, the network gets more dense), the baseline performs more poorly as compared to the proposed approach. 


\section{Conclusion}\label{sec:conclusion}
In this paper we proposed a novel approach for determining a partition of decision variables while decomposing a large-scale optimization problem in a way that improves the performance of distributed optimization methods. We first showed that the partitioning of the decision variables in dual decomposition could be crucial for the empirical performance of a distributed sub-gradient method. Then, we proposed a new method for finding a partition of variables that minimizes the number of constraints being dualized. Our method uses community detection from physics literature to find communities of variables that should be in the same sub-problem in order to achieve the best performance of distributed methods. 

The experimental study using the real application shows that the proposed approach can be used to find a partition for dual decomposition that speeds up the convergence of distributed sub-gradient methods and that the advantage of our approach becomes more significant as each constraint involves more variables and thus, the connectivity among the variables gets stronger. In addition, the proposed methodology can be easily combined with other established techniques that improve the rate of convergence, such as incremental methods \cite{bertsekas2011incremental}, smoothing techniques \cite{nesterov2005smooth,boyd_et_al2011}, adaptive subgradient methods \cite{duchi2011adaptive} among others.  

We highlight here that research in computer science has introduced approaches and algorithms for sub-problem decompositions in a way that communication between computing nodes are minimized \cite{nowak2003,wolfe_et_al2008,knobe1990data,hromkovivc2013communication}; however, there are some key differences between those works and this paper. First, our goal for finding a decomposition is not to minimize communication but to minimize the number of constraints being dualized, in order to conserve as much structure of the original problem as possible while decomposing the optimization problem. Also, each node in the network of this paper represents not a computing node but a decision variable or a group of decision variables (such as a demand location in the transportation problem). While we focused on speeding up the convergence of the distributed method in this paper, the proposed methodology may also be used for minimizing communication between computing nodes, which is a future research topic.

Another potential future research is examining whether the proposed variable partitioning method can be applied to other distributed optimization approaches. Coordinate descent methods \cite{palomar_chiang2006,richtarik_takac2016,wright2015} have gained popularity recently. In coordinate descent methods, variables are partitioned into groups, one of which is chosen to be updated in each iteration. Thus, the approach proposed in this paper can also be used to find a partition for coordinate descent methods; which may also benefit from the community structure of decision variables found by our approach.

\bibliographystyle{plain}
\bibliography{DistributedComputing}

\end{document}